\documentclass[11pt]{article}
\usepackage{amsmath,amsfonts,amssymb,amsthm}
\newtheorem{thm}{Theorem}

\newtheorem{prop}{Proposition}

\def\qed{\nobreak\quad\raise -2pt\hbox{\vrule\vbox to 10pt{\hrule width 6pt
\vfill\hrule}\vrule}\par\vspace{2ex}}
\def\qel{\nobreak\quad\raise -2pt\hbox{\vrule\vbox to 10pt{\hrule width 6pt
\vfill\hrule}\vrule}}

\def\N{\mathbb{ N}}

\def\Z{\mathbb{Z}}

\begin{document}

\begin{center}
{\Large\bf A $q$-Analog of  Dual Sequences with  Applications}
\end{center}

\vskip 2mm \centerline{Sharon J. X. Hou$^a$ and Jiang Zeng$^{b}$}

\begin{center} $^a$Center for Combinatorics, LPMC\\
Nankai University, Tianjin 300071, People's Republic of China\\
{\tt houjx@mail.nankai.edu.cn}
\end{center}

\begin{center}
$^b$Institut Camille Jordan,
Universit\'e Claude Bernard (Lyon I)\\
F-69622 Villeurbanne Cedex, France \\
{\tt zeng@math.univ-lyon1.fr}
\end{center}
{\small \vskip 0.7cm \noindent{\bf Abstract.}   In the present
paper combinatorial identities involving $q$-dual sequences or
polynomials with coefficients that are $q$-dual sequences are
derived. Further, combinatorial identities for $q$-binomial
coefficients(Gaussian coefficients), $q$-Stirling numbers and
$q$-Bernoulli numbers and polynomials are deduced.}

\vskip 0.2mm \noindent{\it Keywords}: $q$-dual sequence,
$q$-binomial coefficients, $q$-Stirling numbers, $q$-Bernoulli
numbers, $q$-Bernoulli polynomials

\vskip 0.2mm \noindent{\bf MR Subject Classifications}: Primary
05A30; Secondary 33D99;

\section{Introduction}
Given a sequence $a_0,a_1,a_2,\ldots, a_n,\ldots$ of elements of a
commutative ring $R$ (for example, the complex numbers,
polynomials or rational functions), one usually describes as
Euler-Seidel matrix associated with $(a_n)$ the double sequence
$(a_n^k)$ ($n\geq 0, k\geq 0$) given by the recurrence~\cite{Du}:
$$
a_n^0=a_n,\quad a_n^k=a_n^{k-1}+a_{n+1}^{k-1}\quad (k\geq 1,n\geq 0).
$$
The sequence $(a_n^0)$  of the first row of the matrix  is the \emph{initial sequence}.
The sequence $(a_0^n)$ of the first column of the matrix is the \emph{final sequence}.
Such a matrix is equivalent to the
  table obtained by computing the finite difference of consecutive terms of $(a_0^n)$ and iterating the procedure.
One passes from the initial sequence to the last one and
conversely through
\begin{equation}\label{eq:euler}
a_0^n=\sum_{i=0}^n{n\choose i}a_i^0\Longleftrightarrow
a_n^0=\sum_{i=0}^n(-1)^{n-i}{n\choose i}a_0^i.
\end{equation}
If one sets  $a_n=(-1)^na_n^0$ and $a_n^*=(-1)^na_0^n$,
then the above relations can be written as
\begin{equation}\label{eq:dual}
a_n^*=\sum_{i=0}^n(-1)^i{n\choose i}a_i\Longleftrightarrow a_n=\sum_{i=0}^n(-1)^i{n\choose i}a_i^*.
\end{equation}
 In \cite{Sun} the sequence $(a_n^*)$ is  called
 the \emph{dual sequence} of $(a_n)$.
It is well-known that if $a_n=(-1)^nB_n$, where $(B_n)=(1, -1/2,
1/6, 0, -1/30, \ldots)$ is the sequence of Bernoulli numbers, then
$a_n^*=a_n$, that is $((-1)^nB_n)$ is \emph{self-dual}.
Generalizing the results of Kaneko~\cite{Ka} and
Momiyama~\cite{Mo} on Bernoulli numbers, Sun~\cite{Sun} has
recently proved some remarkable identities on dual sequences.
 Other generalizations of Kaneko's identity have been obtained by
   Gessel~\cite{Ge} using umbral calculus.

The aim of this paper is to give a $q$-version  of Sun's results in \cite{Sun}.
 In the last two decades there has been an increasing interest
  in  generalizing the classical
  results  with a generic parameter $q$, which is the  so-called phenomenon of "$q$-disease".
 As regards Euler-Seidel matrix Clarke et al.~\cite{CHZ} have given a
 $q$-analog of (1) with application to $q$-enumeration of derangements.

We shall need some standard $q$-notation, which can be found in
Gasper and Rahman's book~\cite{GR}. The $q$-shifted factorial
$(a;q)_n$ is defined by $(a;q)_0=1$ and
$$
(a;q)_n=
(1-a)(1-aq) \cdots (1-aq^{n-1})
$$
if $n$ is a positive integer.
For $k\in \Z$ the $k$-integer $[k]_q$ is defined by $[k]_q= \frac{1-q^k}{1-q}$,
so $[-k]_q=-q^{-k}[k]_q$. For integer $k$,
the $q$-binomial coefficient ${\alpha \brack k}$ is defined by  ${\alpha \brack k}=0$ if $k <0$ and
\[
{\alpha \brack k}=\frac{(1-q^{\alpha})(1-q^{\alpha-1})\cdots (1-q^{\alpha-k+1})}{(q;q)_k}
\]
if $k$ is a positive integer.
Let $(a_n)$ be a sequence of a commutative ring. We call the sequence
$(a_n^*)$ given by
\begin{equation}\label{eq:qdef}
a_n^*=\sum_{i=0}^n{n \brack i}(-1)^ia_iq^{{i \choose 2}}
\end{equation}
the $q$-\emph{dual sequence} of $(a_n)$. By Gauss inversion~\cite[p. 96]{Ai} we get
\begin{equation}\label{eq:qadef}
 a_n=\sum_{r=0}^n{n \brack r}(-1)^ra_r^*q^{{r+1 \choose 2}-nr}.
\end{equation}

We will need the following $q$-analog of binomial formula~\cite[p.
36]{And}:
\begin{equation}\label{qbinomial}
(z;q)_n=\sum_{j=0}^n{n\brack j}(-1)^jz^jq^{j\choose 2},
\end{equation}
and  the $q$-Chu-Vandermonde formula~\cite[p.354]{GR}:
\begin{equation}\label{eq:q-chu-v}
{}_2 \Phi_1\!\left[\begin{matrix}q^{-n}, a\\[5pt]
c \end{matrix};q, q\right]:= \sum _{k=0}^{\infty} \frac
{(q^{-n};q)_k ( a; q)_k}{(c;q)_k}\frac{z^k}{(q;q)_k} ={(c/a;q)_n
\over (c;q)_n}a^n.
\end{equation}

The following is our basic theorem.
\begin{thm} For $k,l\in \N$ the following identities hold true:
\begin{eqnarray}
&&\sum_{j=0}^l {l \brack
j}\frac{(-1)^ja_{k+j+1}^*}{[k+j+1]_q}q^{{j+1 \choose 2}-l(k+j+1)}+
\sum_{j=0}^k {k \brack j}\frac{(-1)^ja_{l+j+1}}{[l+j+1]_q}q^{j+1
\choose 2 }\nonumber\\
&&\hspace{2cm}={a_0 \over [k+l+1]_q{k+l \brack k}} ,
\label{eq:2.1}\\
&& \sum_{j=0}^l{l \brack j}(-1)^ja_{k+j}^*q^{{j+1 \choose
2}-l(k+j)}= \sum_{j=0}^k{k \brack j}(-1)^ja_{l+j}q^{{j
\choose 2}},\label{eq:2.2}\\
&&\sum_{j=0}^{l+1} {l+1 \brack j}(-1)^{j+1}[k+j+1]_q a_{k+j}^* q^{{j
\choose 2}-l(k+j)-k}\nonumber\\
&&\hspace{2cm} = \sum_{j=0}^{k+1} {k+1 \brack j}(-1)^j[l+j+1]_q
{a_{l+j}}q^{j-1 \choose 2}.\label{eq:2.3}
\end{eqnarray}
\end{thm}

The above theorem  is a $q$-analog of
Theorem~2.1 in Sun~\cite{Sun}. Note also that  \eqref{eq:2.2}
was also a $q$-analog of Theorem~7.4 in  Gessel~\cite{Ge}.

The rest of this paper will be organized as follows: we prove Theorem~1 in
Section~2 and present a $q$-analog of Sun's main theorem in
Section~3. In Section~4, we present  some interesting examples as
applications of our Theorems 1 and 2.

\section{Proof of Theorem 1}
Plugging \eqref{eq:qdef} into the first sum of the left-hand side
of \eqref{eq:2.1}, we have
\begin{equation}
LHS =a_0B+\sum_{j=0}^k {k \brack
j}(-1)^j\frac{a_{l+j+1}}{[l+j+1]_q}q^{{j+1 \choose 2
}}+C\label{eq:th5},
\end{equation}
where
\[
B=\sum_{j=0}^l{l \brack j}(-1)^j{q^{{j+1 \choose 2}-l(k+j+1)}\over
[k+j+1]_q},
\]
and
\[
C=\sum_{j=0}^l {l \brack j}(-1)^j{q^{{j+1 \choose
2}-l(k+j+1)}\over [k+j+1]_q}\sum_{i=1}^{k+j+1}{k+j+1 \brack
i}(-1)^ia_iq^{i \choose 2}.
\]
It is known  (see \cite{Zeng}  for further
applications) that
\begin{equation}\label{eq:zeng}
\frac{1}{(x+a_0)(x+a_1) \cdots (x+a_{l})}= \sum_{j=0}^l
\frac{\prod_{i=0 \atop i \neq j}^l(a_i-a_j)^{-1}}{x+a_j}.
\end{equation}
Setting  $x=-q^{-k-1}$ and $a_i=q^i$ ($0\leq i\leq l$) in \eqref{eq:zeng} we obtain
\begin{equation}\label{eq:lagrange}
\sum_{j=0}^l(-1)^j\frac{q^{{j+1 \choose
2}-l(k+j+1)}}{(q;q)_j(q;q)_{l-j}(1-q^{k+j+1})}=\frac{1}{(q^{k+1};
q)_{l+1}}.
\end{equation}
It follows that
$$B=\frac{(1-q)(q;q)_l}{(q^{k+1};q)_{l+1}}=\frac{1}{
[k+l+1]_q{k+l \brack k}}.
$$
Exchanging the order of summation we can rewrite $C$ as follows:
\begin{eqnarray*}
C&=&\sum_{i=1}^{k+l+1}(-1)^i {a_i \over [i]_q}q^{i \choose
2}\sum_{j=i-k-1}^l(-1)^j {l \brack j}{k+j \brack i-1}q^{{j+1
\choose 2}-lj-(k+1)l}\\
&=& \sum_{i=1}^k(-1)^i {a_i \over [i]_q}q^{i \choose 2 } {k
\brack i-1}{}_2 \Phi_1\!\left[\begin{matrix}q^{-l}, q^{k+1}\\[5pt]
q^{k-i+2}\end{matrix};q, q\right]q^{-(k+1)l}.
\end{eqnarray*}
Applying the $q$-Chu-Vandermonde
formula \eqref{eq:q-chu-v} we obtain
\begin{eqnarray*}
C&=&  \sum_{i=1}^k (-1)^i {a_i \over [i]_q}q^{i \choose 2 } {k
\brack i-1}\frac{(q^{-i+1};q)_l}{(q^{k-i+2};q)_l}\\
&=& \sum_{i=1}^k(-1)^{i+l}{k \brack i-l-1}{a_i \over [i]_q}q^{{i
\choose 2}+{l+1 \choose 2}-il}\\
&=& \sum_{j=0}^k {k \brack
j}(-1)^{j+1}\frac{a_{l+j+1}}{[l+j+1]_q}q^{{j+1 \choose 2 }}.
\end{eqnarray*}
Substituting the values of $B$ and  $C$ into \eqref{eq:th5} yields \eqref{eq:2.1}.

To derive \eqref{eq:2.2} and \eqref{eq:2.3} from
 \eqref{eq:2.1} we
define the linear operator $\delta_q$ by
$$
\delta_q(a_n)=-q^{1-n}[n]_q a_{n-1}\qquad\textrm{for}\quad n\geq 0.
$$
Then  $\delta_q(a_n^*)=[n]_qa_{n-1}^*$.
 Indeed,
\begin{eqnarray*}
\delta_q(a_n^*)&=&\sum_{i=0}^n{n \brack i}(-1)^i\delta_q(a_i)q^{i \choose 2}
 = \sum_{i=0}^n{n \brack i}(-1)^{i+1}q^{1-i}[i]_qa_{i-1}q^{i \choose 2}\\
&=&[n]_q\sum_{i=0}^n(-1)^i {n-1 \brack i-1}(-1)^{i-1}a_{i-1}q^{i-1 \choose 2} = [n]_qa_{n-1}^*.
\end{eqnarray*}
Now, applying $\delta_q$ to \eqref{eq:2.1} yields \eqref{eq:2.2}.
Furthermore, replacing $k$ by $k+1$ and $l$ by $l+1$ in \eqref{eq:2.2}
then applying $\delta_q$ on both sides yields \eqref{eq:2.3}.

\noindent{\bf Remark}:
We can also prove \eqref{eq:2.2} and \eqref{eq:2.3} directly by using
the $q$-Chu-Vandermonde formula.

\section{A $q$-analog of Sun's main theorem}
In this section, we assume that $x$, $y$  and $z$ are commuting
indeterminates. Define $[x,y]^n$ by $[x,y]^0=1$ and
$$
[x, y]^n=\sum_{i=0}^n{n\brack i} x^iy^{n-i}
$$
for positive integer $n$. So $[x,y]^n=(x+y)^n$ when $q=1$.
Similarly
$$
[x,y,z]^n=[x,[y,z]]^n=\sum_{i=0}^n{n\brack i}
x^i[y,z]^{n-i}=\sum_{i,j,k\geq 0\atop i+j+k=n
}\frac{[n]_q!}{[i]_q![j]_q![k]_q!}x^iy^jz^k,
$$
and hence $[x,y,z]^n$ is a symmetric polynomial of $x,y,z$ and
$[x,y,z]^n=(x+y+z)^n$ when $q=1$.

 Like the definition of Bernoulli polynomials, we introduce
$$
A_n(x)=\sum_{i=0}^n(-1)^i{n \brack i}a_iq^{i \choose 2}x^{n-i}
\quad \text{and} \quad A_n^*(x)=\sum_{i =0}^n(-1)^i{n \brack
i}a_i^*x^{n-i}.
$$
The following is our $q$-analog of the main theorem of
Sun~\cite[Th. 1.1]{Sun}.
\begin{thm}
Let $k, l \in \N$, then
\begin{eqnarray}
&&(-1)^l\sum_{j=0}^l{l \brack j}x^{l-j}{A_{k+j+1}^*(z) \over
[k+j+1]_q}q^{-kj-{k+1 \choose 2}} \nonumber\\
 &&\hspace{1cm} +(-1)^k\sum_{j=0}^k{k \brack
j}x^{k-j}{A_{l+j+1}([1,-z,-x]) \over [l+j+1]_q}q^{{j+1 \choose
2}-k(l+j+1)}
=\frac{a_0(-x)^{k+l+1}}{[k+l+1]{k+l \brack k}}.\label{eq:qsun1}\\
&& (-1)^l\sum_{j=0}^l{l \brack j}x^{l-j}A_{k+j}^*(z)q^{k(l-j)}
=(-1)^k\sum_{j=0}^k{k \brack j}x^{k-j}A_{l+j}([1,-z,-x])
q^{k-j \choose 2}.\label{eq:qsun2} \\
&&(-1)^{l+1}\sum_{j=0}^{l+1}{l+1 \brack j}x^{l+1-j}[k+j+1]_q
A_{k+j}^*(z) q^{(k+1)(l-j)+1}\nonumber \\
&& \hspace{1cm}=(-1)^k\sum_{j=0}^{k+1}{k+1 \brack
j}x^{k+1-j}[l+j+1]_q A_{l+j}([1,-z,-x]) q^{{k-j \choose 2}-j}.
\label{eq:qsun3}
\end{eqnarray}
\end{thm}

\begin{proof}
We derive from \eqref{eq:qadef} and \eqref{qbinomial} that
\begin{eqnarray}\nonumber
A_n([1,-x])&=&\sum_{i=0}^n{n \brack i}(-1)^ia_i[1,-x]^{n-i}q^{i
\choose 2}\\ \nonumber &=&\sum_{i,j,s \geq 0}^n{n \brack i}{i
\brack j}(-1)^{i-j}a_j^*{n-i \brack s}(-x)^sq^{{i-j \choose 2}}\\
\nonumber &=&\sum_{j = 0}^n{n \brack j}a_j^*\sum_{i,s \geq 0}{n-j
\brack s}(-1)^sx^s {n-j-s \brack i-j}(-1)^{i-j}q^{{i-j \choose
2}}\\ \nonumber  &=&(-1)^n\sum_{j=0}^n{n \brack
j}(-1)^ja_j^*x^{n-j}\\ \label{eq:qax} &=&(-1)^nA_n^*(x),
\end{eqnarray}
and
\begin{eqnarray}\nonumber
A_n([1,-z,-x])&=&\sum_{i=0}^n{n \brack
i}(-1)^ia_i[1,-z,-x]^{n-i}q^{{i \choose 2}}\\ \nonumber
&=&\sum_{i, j \geq 0}{n \brack i}{n-i \brack
j}(-1)^{i+j}x^ja_i[1,-z]^{n-i-j}q^{i \choose 2}\\ \nonumber
&=&\sum_{i, j \geq 0}{n \brack
j}(-1)^jx^{j}{n-j \brack i}(-1)^ia_i[1,-z]^{n-i-j}q^{i \choose 2}\\
\nonumber &=&\sum_{j=0}^n{n \brack
j}(-1)^jx^jA_{n-j}([1,-z])\\\label{eq:qinver} &=&
(-1)^n\sum_{j=0}^n{n \brack j}x^jA_{n-j}^*(z).
\end{eqnarray}
Denote the first sum of the left-hand side in \eqref{eq:qsun1} by
$\mathcal{C}$.
 Applying \eqref{eq:qax} and \eqref{eq:qinver},
the left-hand side of \eqref{eq:qsun1} is equal to
\begin{eqnarray}\nonumber
&& (-1)^k\sum_{j=0}^k{k \brack j}{x^{k-j} \over [l+j+1]_q}q^{{j+1
\choose 2}-k(l+j+1)}
\sum_{i=0}^{l+j+1}{l+j+1 \brack i}A_i([1,-z])(-x)^{l+j+1-i} +\mathcal{C}\\
&=&a_0(-x)^{k+l+1}\mathcal{B}+S+\mathcal{C},
\end{eqnarray}
where
\[
\mathcal{B}=\sum_{j=0}^k{k \brack j}(-1)^j{q^{{j+1 \choose
2}-k(l+j+1)} \over [l+j+1]_q},
\]
and
$$
S=(-1)^{k+l+1}\sum_{j=0}^k{k \brack j}(-1)^jq^{{j+1 \choose
2}-k(l+j+1)}\sum_{i=1}^{l+j}{l+j \brack
i-1}x^{k+l+1-i}{A_i^*(z)\over [i]_q}.
$$
Exchanging  $k$ and $l$ in \eqref{eq:lagrange} yields
$$
 \mathcal{B}={1 \over [k+l+1]_q {k+l \brack k}}.
$$
Now, we show that $S=-\mathcal{C}$. Exchanging the order of
summation we have
\begin{eqnarray*}
S&=&(-1)^{k+l+1}\sum_{i=1}^l{l \brack i-1}{A_i^*(z) \over
[i]_q}x^{k+1+l-i}
{}_2 \Phi_1\!\left[\begin{matrix} q^{-k}, q^{l+1}\\[5pt]
q^{l-i+2}\end{matrix};q, q\right]q^{-k(l+1)}\\
&=&(-1)^{k+l+1}\sum_{i=1}^l{l \brack i-1} {A_i^*(z) \over
[i]_q}x^{k+l-i} {(q^{-i+1}; q)_k \over (q^{l-i+2};q)_k}
\ \ \text{(by $q$-Chu-Vandemonde)}\\
&=&(-1)^{l+1}\sum_{i=1}^l{l \brack i-k-1}{A_i^*(z) \over [i]_q}x^{k+1+l-i}q^{-ik+{k+1 \choose 2}}\\
&=&(-1)^{l+1}\sum_{j=0}^l{l \brack j}{A_{k+j+1}^*(z) \over
[k+j+1]_q}x^{l-j}q^{-jk-{k+1 \choose 2}}=-\mathcal{C}.
\end{eqnarray*}
Next, the right-hand side of \eqref{eq:qsun2} is equal to
\begin{eqnarray*}
& &(-1)^k\sum_{j=0}^k{k \brack j}x^{k-j}\sum_{i=0}^{l+j}{l+j
\brack
i}(-x)^{l+j-i}A_i([1,-z])q^{k-j \choose 2}\\
&=&(-1)^{k+l}\sum_{i=0}^l{l \brack i}x^{k+l-i}A_i^*(z)\ {}_2 \Phi_1\!\left[\begin{matrix}q^{-k}, q^{l+1}\\[5pt]
q^{l-i+1}\end{matrix};q, q\right]q^{k \choose 2}\\
&=&(-1)^{k+l}\sum_{i=0}^l{l \brack
i}x^{k+l-i}A_i^*(z)\frac{(q^{-i};q)_k}{(q^{l-i+1};q)_k}q^{{k \choose 2}+k(l+1)}\\
&=& (-1)^{l}\sum_{i=0}^l{l \brack i-k}x^{k+l-i}A_i^*(z)q^{-ik+k^2+kl}\\
&=&(-1)^{l}\sum_{j=0}^l{l \brack j}x^{l-j}A_{k+j}^*(z)q^{k(l-j)},
\end{eqnarray*}
which is exactly the left-hand side of \eqref{eq:qsun2}.

Finally, exchanging the order of summation, the right-hand side of
\eqref{eq:qsun3} can be written as
\begin{eqnarray*}
R&= &(-1)^k\sum_{j=0}^{k+1}{k+1 \brack j}x^{k+1-j}[l+j+1]_q
\sum_{i=0}^{l+j}{l+j \brack
i}(-x)^{l+j-i}A_i([1,-z])q^{{k-j \choose 2}-j}\\
&=&(-1)^{k+l}\sum_{i,j \geq 0}{k+1 \brack j}x^{k+l+1-i}[i+1]_q
(-1)^j{l+j+1 \brack i+1}A_i^*(z)q^{{k-j \choose 2}-j}\\
&=&(-1)^{k+l}\sum_{i=0}^{l+1}{l+1 \brack i}x^{k+l+1-i}[i+1]_qA_i^*(z) \ {}_2 \Phi_1\!\left[\begin{matrix}q^{-k-1}, q^{l+2}\\[5pt]
q^{l-i+1}\end{matrix};q, q\right]q^{k \choose 2}.
\end{eqnarray*}
By $q$-Chu-Vandermonde formula we have
\begin{eqnarray*}
R&=&(-1)^{k+l}\sum_{i=0}^{l+1}{l+1 \brack
i}x^{k+l+1-i}[i+1]_qA_i^*(z)\frac{(q^{-i-1};q)_{k+1}}{(q^{l-i+1};q)_{k+1}}q^{{k \choose 2}+(k+1)(l+2)}\\
&=& (-1)^{l+1}\sum_{i=0}^{l+1}{l+1 \brack i-k}x^{k+l+1-i}[i+1]_qA_i^*(z)q^{(k+1)(l+2-i)+k^2-k-1}\\
&=&(-1)^{l+1}\sum_{j=0}^{l+1}{l+1 \brack
j}x^{l+1-j}[k+j+1]_qA_{k+j}^*(z)q^{(k+1)(l-j)+1},
\end{eqnarray*}
which is exactly the left-hand side of \eqref{eq:qsun3}.
\end{proof}

\noindent {\bf Remark.} When $q=1$, Theorems 1 and  2, which correspond to
Theorems 2.2 and 1.1 of Sun\cite{Sun}, are actually equivalent. Indeed, in such case,
we have
\begin{equation}\label{eq:lem}
(-1)^nA_n^*(1-x)=A_n(x),
\end{equation}
which can be verified as follows:
\begin{eqnarray*}
\sum_{i=0}^n{n\choose i}a_i^*(x-1)^{n-i}&=& \sum_{i=0}^n{n\choose
i}(-1)^{n-i}\sum_{j=0}^i{i\choose
j}(-1)^ja_j(1-x)^{n-i}\\
&=&\sum_{j=0}^n{n\choose j}(-1)^ja_j\sum_{i=j}^n{n-j\choose
i-j}(x-1)^{n-i}\\
&=&\sum_{j=0}^n{n\choose j}(-1)^ja_jx^{n-j}.
\end{eqnarray*}
Now, taking  $a_n=(-1)^{l+k+n}x^{k+l-n}A_n(y)$ with $q=1$,
\begin{eqnarray*}
a_n^*&=& \sum_{i=0}^n{n\choose i}(-1)^{l+k}x^{k+l-i}A_i(y)\\
&=&\sum_{i=0}^n{n\choose
i}(-1)^{l+k}x^{k+l-i}\sum_{j=0}^i{i\choose j}(-1)^ja_jy^{i-j}\\
&=&(-1)^{l+k+n}x^{k+l-n}\sum_{j=0}^n{n\choose
j}a_j(-1)^{n-j}\sum_{i=j}^n {n-j\choose i-j}x^{n-i}y^{i-j}\\
&=&(-1)^{l+k}x^{k+l-n}A_n(x+y).
\end{eqnarray*}
It follows from  \eqref{eq:lem} that
$$
a_n^*=(-1)^{k+l+n}x^{k+l-n}A_n^*(1-x-y).
$$
Substituting the above values of $a_n$ and $a_n^*$ in Theorem~1 we
obtain Theorem~2.
Conversely, it is easy to see that
 Theorem~1 is a special case of Theorem~2 because
$$
A_n(0)=(-1)^na_n, \qquad A_n(1)=a_n^*.
$$
Hence we have proved that Theorems~1 and 2 are actually
equivalent when $q=1$.
\section{Some applications}
In this section we derive some examples from our main theorem, most of them are
 $q$-analogs of results in Sun~\cite{Sun}.
\subsection*{Example 1}
For any fixed integer $i\geq 0$ let $a_n=(-1)^n{n \brack
i}t^{n-i}q^{{i \choose 2}}$, then it follows from \eqref{eq:qdef}
and \eqref{qbinomial} that
\begin{eqnarray*}
a_n^*&=&\sum_{k=0}^n{n\brack k}{k\brack i}t^{k-i}q^{{i\choose 2}+{k\choose 2}}\\
&=& {n\brack i}q^{i^2-i}\sum_{k=i}^n{n-i\brack k-i}(tq^i)^{k-i}q^{k-i\choose 2}\\
&=&{n \brack i}(-tq^i;q)_{n-i}q^{i^2-i}.
\end{eqnarray*}
Substituting the above values in \eqref{eq:2.2} of  Theorem 1 yields
\begin{equation}\label{qprop}
\sum_{j=0}^l{l\brack j}{k+j\brack
i}(-1)^{l-j}(-q^it;q)_{k+j-i}q^{j(j+1)/2-lj+{i\choose 2}}=
\sum_{j=0}^k{k\brack j}{l+j\brack i}t^{l+j-i}q^{kl+j(j-1)/2}.
\end{equation}
For variations of methods, we will give two more proofs of
\eqref{qprop}. Note that when $q=1$ Eq.~\eqref{qprop} reduces to a
crucial result of Sun \cite[Lemma 3.1]{Sun}, which was proved  by
using derivative operator.

We  first $q$-generalize Sun's proof by
using $q$-derivative operator.
  For any polynomial $f(t)$ in $t$, let $D_q$ be the $q$-derivative operator with respect to $t$:
$$
D_qf(t)=\frac{f(tq)-f(t)}{(q-1)t}.
$$
Clearly we have
$$
D_qt^n=\frac{q^n-1}{q-1}t^{n-1},\quad D_q((-t;q)_n)=[n]_q(-qt;q)_{n-1}.
$$
For integer $i\geq 0$ define $[i]_q!=\prod_{j=0}^i[j]_q$, then
\begin{eqnarray}
D_q^i(t^n)&=&[i]_q!{n\brack i} t^{n-i},\label{eq:D1}\\
D_q^i((-t;q)_n)&=&q^{i(i-1)/2}[i]_q!{n\brack i}(-q^it;q)_{n-i}.\label{eq:D2}
\end{eqnarray}
By Gauss inversion, the $q$-binomial formula~\eqref{qbinomial} is
equivalent to
$$
z^n=\sum_{j=0}^n(-1)^jq^{{j+1\choose 2}-nj}{n\brack j}(z;q)_j.
$$
Substituting $z$ by $-tq^k$ we get
\begin{equation}\label{eq:inverseqbin}
(tq^k)^n=\sum_{j=0}^n{n\brack
j}(-1)^{n-j}q^{j(j+1)/2-nj}(-tq^k;q)_j.
\end{equation}
Now, using the $q$-derivative operator and
\eqref{eq:D1}-\eqref{eq:inverseqbin}, we can write the difference
of the two sides of \eqref{qprop} as follows:
\begin{eqnarray*}
&&\frac{1}{[i]_q!}D^i_q\left((-t;q)_k\sum_{j=0}^l{l\brack
j}(-1)^{l-j}q^{j(j+1)/2-lj}(-tq^k;q)_j-(tq^k)^l\sum_{j=0}^k{k\brack j}q^{j(j-1)/2}t^j\right)\\
&=&\frac{1}{[i]_q!}D^i_q\left((-t;q)_k(tq^k)^l-(tq^k)^l(-t;q)_k\right),
\end{eqnarray*}
which is clearly equal to 0.

Our second proof of \eqref{qprop} uses  the machinery of basic
hypergeometric functions. Rewriting \eqref{qprop} in terms of basic
hypergeometric functions, we have
\begin{eqnarray}\label{eq:basic}
&&{k\brack i}(-1)^l(-q^it;q)_{k-i}q^{i\choose 2}{}_3\phi_2\left[
\begin{array}{c}q^{-l},\; q^{k+1},\;-tq^k\\
q^{k-i+1},\; 0\end{array};q,q\right]\nonumber\\
&&\hspace{2cm}={l\brack
i}t^{l-i}q^{kl}{}_2\phi_1\left[
\begin{array}{c}
q^{-k},\;q^{l+1}\\
q^{l-i+1}\end{array};q, -tq^k \right].
\end{eqnarray}
A standard proof of \eqref{eq:basic} goes then as follows:
\begin{eqnarray*}
&&{k \brack i}(-1)^l(-tq^i;q)_{k-i}q^{{i \choose
2}}(-tq^{k-i})^l
{}_3\phi_2\!\left[\begin{matrix}q^{-l}, q^{-i},(-tq^{i-1})^{-1}\\[5pt]
q^{k-i+1}, 0\end{matrix}; q, q\right] \\
&&\hspace{3cm}\text{(by  ~\cite[p.241(III.11)]{GR})}\\
& = & {k \brack i} t^lq^{(k+i)l}q^{i \choose
2}(-tq^i;q)_{k-i}\frac{(q^{l-i+1};q)_i}{(q^{-k};q)_i}(tq^{k+l})^{-i}
{}_2\phi_1\!\left[\begin{matrix}q^{-i}, q^{k+l+1-i}\\[5pt]
q^{l-i+1}\end{matrix}; q, -tq^i\right]\\
&&\hspace{3cm} \text{(by  ~\cite[p.241(III.6)]{GR})}\\
& = & {l \brack i}q^{kl}t^{l-i}(-tq^i;q)_{k-i}
\frac{(-tq^k;q)_\infty}{(-tq^i;q)_\infty}{}_2\phi_1\!
\left[\begin{matrix}q^{-k}, q^{l-1}\\[5pt]
q^{l-i+1}\end{matrix}; q, -tq^k\right] \\
&&\hspace{3cm}\text{(by  ~\cite[p.241(III.3)]{GR})}
\end{eqnarray*}
which is equal to the right-hand side of \eqref{eq:basic}.
\subsection*{Example 2}
Let $a_n={x+n \brack m}q^{-mn}$ for $n \in \N$. By the notation
\eqref{eq:qdef}, we have
\[
a_n^*=\sum_{i=0}^n{n \brack i}(-1)^i{x+i \brack m}q^{{i \choose 2}-im}=\left\{
\begin{array}{lr}
(-1)^n{x \brack m-n}q^{-mn+{n \choose 2}} & \quad \text{if}\  m \geq n,\\
0 &\quad \text{otherwise}.\\
\end{array}
\right.
\]
 Theorem~1 implies that
\begin{eqnarray*}
& &\sum_{j=0}^k{k \brack j}{(-1)^j \over [l+j+1]_q}{x+l+j+1 \brack
m}q^{-m(l+j+1)+{j \choose 2}} \\
&=& (-1)^{k}\sum_{k \leq j \leq m}{1 \over [j]_q}{l \brack j-k-1}{x
\brack m-j}q^{{j-k \choose 2}+{j \choose 2}-l(j-1)-mj} + \frac{{x
\brack m}}{[k+l+1]_q{k+l \brack k}}.
\end{eqnarray*}

\subsection*{Example 3}
Let $c_n={y \brack n}/{x \brack n}$ for $n \in \N$. Then
$c_n^*={x-y \brack n}q^{ny}/{x \brack n}$. In fact,
\[
{n \brack k}=(-1)^k{-n+k-1 \brack k}q^{nk-{k \choose 2}}.
\]
\begin{eqnarray*}
{x \brack n}c_n^*&=&\sum_{i=0}^n{x-i \brack n-i}(-1)^i{y \brack i}q^{{i \choose 2}}
=(-1)^n\sum_{i=0}^n{x-n+1 \brack n-i}{y \brack i}q^{(x-i)(n-i)-{n-i \choose 2}+{i \choose 2}}\\
&=&(-1)^nq^{xn-{n \choose 2}}{n-x-1+y \brack n}={x-y \brack n}q^{ny}.
\end{eqnarray*}
By the identities in Theorem~1, we obtain
\begin{eqnarray*}
\sum_{j=0}^k{k \brack j}\frac{(-1)^j{y \brack l+j+1}}{{x-1 \brack l+j}}q^{j+1 \choose 2}
+\sum_{j=0}^l\frac{(-1)^j{x-y \brack k+j+1}}{{x-1 \brack k+j}}q^{(k+j+1)y+{j+1 \choose 2}-l(k+j+1)}
=\frac{[x]_q}{[k+l+1]{k+l \brack k}}
\end{eqnarray*}
and
\begin{eqnarray*}
\sum_{j=0}^k{k \brack j}(-1)^j\frac{{y \brack l+j}}{{x \brack l+j}}q^{j \choose 2}=
\sum_{j=0}^l(-1)^j\frac{{x-y \brack k+j}}{{x \brack k+j}}q^{(k+j)y+{j+1 \choose 2}-l(k+j)}.
\end{eqnarray*}

\subsection*{Example 4}
 Carlitz~\cite[(3.1)]{Car} defined the \emph{$q$-Stirling numbers of the second kind} ${m \brace
n}_q$ by
\[
[n]_q^m=\sum_{i=0}^m {m \brace i}_q [i]_q!{n\brack i}q^{i \choose 2}.
\]
By Gauss inversion we get
\begin{eqnarray*}
{m \brace n}_q & = & \frac{q^{-{n \choose
2}}}{[n]_q!}\sum_{i=0}^n(-1)^iq^{i
\choose 2}{n \brack i}[n-i]_q^m \\
&=& {1 \over [n]_q!}\sum_{i=0}^n(-1)^{n-i}q^{{i+1 \choose 2}-ni}{n
\brack i}[i]_q^m.
\end{eqnarray*}
So we have the following $q$-dual sequences:
$$
a_n =(-1)^n[n]_q!{m \brace n}_q,\quad
a_n^*=[n]_q^m.
$$
Substituting these values in  Theorem~1 yields corresponding identities. For example,
 applying
\eqref{eq:2.2} we obtain
\begin{eqnarray*}
{1 \over [l]_q!}\sum_{j=0}^l(-1)^{l-j}q^{{j+1 \choose 2}-l(k+j)}{l
\brack j} [k+j]_{q}^m  = \sum_{j=0}^k q^{j \choose 2}{l+j \brack
j}{[k]_q! \over [k-j]_q!}{m \brace l+j}_q .
\end{eqnarray*}
The left-hand side of the above identity is called a \emph{non-central $q$-Stirling number of
the second kind}, with non-centrality parameter $k$, by Charalambides~\cite{Char}. This number was
first  discussed by Carlitz~\cite[(3.8)]{Car} and recently by
Charalambides~\cite[(3.5)]{Char}. Note that for $k=0$ these numbers
reduce to the usual $q$-Stirling numbers of the second kind, while
for $k \neq 0$ the above identity connects the non-central to the
usual $q$-Stirling numbers of the second kind.

\subsection*{Example 5}
Taking
$e(t)=\sum_{n=0}^\infty\frac{x^n}{[n]_q!}$
as a q-analog of the exponential function $e^x$,
Al-Salam~\cite[2.1]{Salam} defined a $q$-analog of Bernoulli numbers $B_n$ by
\[
\frac{1}{e(t)-1}=\sum_{n=0}^\infty\frac{t^n}{[n]_q!}B_n.
\]
These $q$-Bernoulli numbers $B_n$ satisfy the following recurrence
relation (see \cite[4.3]{Salam}):
\[
[1, B]^n =\left\{
\begin{array}{cc}
B_n & n > 1,\\
1+B_1& n=1.\\
\end{array}
 \right.
 \]
Now, if  $a_0^*=B_0=1$ and for $n \geq  1$, $a_n^*=B_n=\sum_{i=0}^n{n \brack
i}B_i$, then $a_0=1$ and for $n\geq 1$
\begin{eqnarray*}
a_n&=&\sum_{i=0}^n {n \brack i}(-1)^iq^{{i+1 \choose
2}-in}\sum_{j=0}^i{i \brack j}B_j\\
&=&\sum_{j=0}^n{n \brack
j}B_j\sum_{i\geq j}{n-j \brack i-j}(-1)^iq^{{i+1 \choose
2}-ni}\\
&=&(-1)^nB_nq^{-{n \choose 2}}.
\end{eqnarray*}
Theorem~1 infers then the following identities, which are
$q$-analogs of the identities of Kaneko~\cite{Ka} and
Momiyama~\cite{Mo} on Bernoulli numbers.
\begin{prop}  For $k, l \in \N$,
\begin{eqnarray*}
& &\sum_{j=0}^k {l \brack j}\frac{(-1)^jB_{k+j+1}}{[k+j+1]_q}q^{{j+1 \choose 2}-l(k+j+1)}
+\sum_{j=0}^l{k \brack j}\frac{(-1)^{l+1}B_{l+j+1}}{[l+j+1]_q}q^{-{l+1 \choose 2}-lj}
=\frac{1}{[k+l+1]_q {k+l \brack k}} ,\\
&& \sum_{j=0}^l{l \brack j}(-1)^jB_{k+j}q^{{l-j \choose 2}}=\sum_{j=0}^k{k \brack j}(-1)^l B_{l+j}q^{l(k-j)},\\
&& \sum_{j=0}^{l+1}{l+1 \brack j}(-1)^{j+1}[k+j+1]_q B_{k+j}q^{{l-j \choose 2}-j}
 =\sum_{j=0}^{k+1} {k+1 \brack j}(-1)^l[l+j+1]_q B_{l+j}q^{(k-j)(l+1)+1}.
\end{eqnarray*}
\end{prop}

\subsection*{Example 6}
Al-Salam~\cite{Salam} also defined the $q$-Bernoulli polynomials
$B_n(x)=\sum_{k=0}^n{n \brack k}B_kx^{n-k}$. By Example~5,
if  $a_n=(-1)^n B_nq^{-{n \choose 2}}$ then  $a_n^*=B_n$. Therefore
$$
A_n(x)=\sum_{i=0}^n{n \brack i}B_ix^{n-i}=B_n(x) \quad \text{and}
\quad A_n^*(x)=\sum_{i=0}^n{n \brack i}(-1)^iB_ix^{n-i}=B_n^*(x).
$$
If we replace $A_n(x)$ and $A_n^*(x)$ by $B_n(x)$ and $B_n^*(x)$,
respectively,
 we get the following result.
\begin{prop}
For $k, l \in \N$,
\begin{eqnarray*}
&&(-1)^l\sum_{j=0}^l{l \brack j}x^{l-j}{B_{k+j+1}^*(z) \over
[k+j+1]_q}q^{-kj-{k+1 \choose
2}}\\
&& \hspace{1cm}+(-1)^k\sum_{j=0}^k{k \brack
j}x^{k-j}{B_{l+j+1}([1,-x,-z]) \over [l+j+1]_q}q^{{j+1 \choose
2}-k(l+j+1)}=\frac{a_0(-x)^{k+l+1}}
{[k+l+1]{k+l \brack k}},\\
&&(-1)^l\sum_{j=0}^l{l \brack j}x^{l-j}B_{k+j}^*(z)q^{k(l-j)}
=(-1)^k\sum_{j=0}^k{k \brack j}x^{k-j}B_{l+j}([1,-x,-z])q^{k-j
\choose
2},\\
&&(-1)^{l+1}\sum_{j=0}^{l+1}{l+1 \brack
j}x^{l+1-j}[k+j+1]_q B_{k+j}^*(z) q^{(k+1)(l-j)+1}\\
&& \hspace{3cm}=(-1)^k\sum_{j=0}^{k+1}{k+1 \brack
j}x^{k+1-j}[l+j+1]_q B_{l+j}([1,-x,-z])q^{{k-j \choose
2}-j}\nonumber.
\end{eqnarray*}
\end{prop}
\begin{remark}
It is easy to see that $B_n(0)=B_n$ and $B_n(0)^*=(-1)^nB_n$. Hence
Proposition 1 can be derived from Proposition 2 by taking $x=1$
and $z=0$.
\end{remark}

\section*{Acknowledgements} This work was done under the auspices of
the National Science Foundation of China. The second author
thanks Sun Zhi-Wei for asking  a $q$-question about  the results
in \cite{Sun}, and was also  supported by EC's IHRP Programme, within Research
Training Network ``Algebraic Combinatorics in Europe'', grant HPRN-CT-2001-00272.


\end{document}